%1/28/98

\font\tenbb=msbm10
\font\sevenbb=msbm7
\font\fivebb=msbm5
\newfam\bbfam
\textfont\bbfam=\tenbb
\scriptfont\bbfam=\sevenbb
\scriptscriptfont\bbfam=\fivebb
\def\bb{\fam\bbfam\tenbb}
\def\frac#1#2{{{#1}\over{#2}}}

\input epsf
\magnification = 1200
\centerline {{\bf A SIMPLY CONNECTED NUMERICAL GODEAUX SURFACE}}
\centerline{{\bf WITH AMPLE CANONICAL CLASS}}
\vglue .4in
\centerline{ I. DOLGACHEV \footnote{*}{{\sevenrm
Research supported in part by a NSF grant}.} and C. WERNER}
\vglue .5in

\noindent
\centerline{\bf 1. Introduction} 
\medskip
This paper is a follow-up of a recent paper of PierCarlo
Craighero and Remo Gattazzo {\bf [CG]} in which a new construction of a numerical Godeaux
surface was given. Recall  that a numerical Godeaux surface is a minimal surface $V$ of
general type with zero geometric genus 
$p_g$ and $K_V^2$ = 1. The surface is obtained as a minimal resolution of a quintic surface in 
${\bb P}^3$ with four simple elliptic singularities of type $z^2+x^3+y^6 = 0$. 
In this paper we show that the surface is simply connected and we also find its birational 
model as a double plane branched along an irreducible curve 
of degree 10 with five singular points of type $(3,3)$ and one quadruple ordinary 
singular point. Previous examples of such a curve ({\bf [R]},{\bf [St]}) lead to 
classical Godeaux surfaces (with fundamental group of order 5). A reducible curve of this
kind was constructed earlier by F. Oort and C. Peters
${\bf [OP]}$ and the second  author {\bf [W]}. 
The corresponding double cover is 
birationally equivalent to a numerical Godeaux surface with fundamental 
group ${\bb Z}/4{\bb Z}$ and ${\bb Z}/2{\bb Z}$, respectively. 

Recall that 
$\{1\},{\bb Z}/2{\bb Z},{\bb Z}/3{\bb Z},{\bb Z}/4{\bb Z},{\bb Z}/5{\bb Z}$ are the only 
possible values for the torsion group
of a numerical Godeaux surface. The first (and the only) construction of  
simply connected surfaces of general type with $p_g = 0$  was given
by  R. Barlow {\bf [B]}. Although Barlow's surfaces are numerical Godeaux surfaces they are
not isomorphic to the surfaces of Craighero-Gattazzo. 
In particular the  Craighero-Gattazzo surfaces are 
the first examples of simply connected surfaces with $p_g = 0$ and ample canonical class.
We prove this by showing that the 
surfaces contain no smooth rational curves with self-intersection $(-2)$; however
each Barlow surface contains four such curves.   
In {\bf [CL]},  F. Catanese and C. LeBrun
prove that the Barlow surface can be deformed to a surface with ample canonical
class.
It is  possible that Craighero-Gattazzo surfaces are deformation equivalent to
Barlow surfaces.    

Recently F. Catanese has announced a new construction of  a family of simply connected
numerical Godeaux surfaces. We do not know whether our surfaces belong to his family.

Finally the first author would like to thank E. Stagnaro for bringing to his attention the
paper of Craighero and Gattazzo, and also M. Reid for useful discussions. Both authors
thank J. Keum for his careful reading of the manuscript and many constructive comments. 

\bigskip
\centerline{{\bf 2. The Quintics}}
\smallskip 
Consider the following projective automorphism of ${\bb
P}^3$:
$$\sigma: (X,Y,Z,T) \to (T,X,Y,Z).$$
Its order is four and its fixed points in ${\bb P}^3$ are
$$P_0 = (1,1,1,1), P_1 = (1,i,-1,-i),  P_2 = (1,-i,-1,i), Q_0 = (1,-1,1,-1).$$
The set of fixed points of 
$$\sigma^2: (X,Y,Z,T) \to (Z,T,X,Y)$$ is equal to the union of two lines
$$ r = \{X+Z = Y+T=0 \} = <P_1,P_2>,\ \ r':\{X-Z = Y-T = 0\} = <P_0,Q_0>.$$
We are looking for a quintic surface which is invariant with respect to $\sigma$ and 
has simple elliptic singularities of degree 1 locally isomorphic to 
$z^2+x^3+y^6 = 0$ (tacnodal points) at the reference points 
$$a_1 = (1,0,0,0),\quad a_2 = (0,1,0,0),\quad a_3 = (0,0,1,0),\quad a_4 = (0,0,0,1).$$
It is easy to check that any homogeneous polynomial $F_5(X,Y,Z,T)$ 
of degree 5 which is invariant with respect to $\sigma$ and has degenerate 
critical points of multiplicity 2 at 
the points $a_i$ must 
look like

\noindent$
\displaystyle{
\; (aT+bZ+Y )^{2}{X}^{3}+ (aX+bT+Z )^{2}{Y}^{3}+
 (aY+bX+ T )^{2}{Z}^{3}+ (aZ+bY+X )^{2}{T}^{3}}+$

\noindent
${
 ({X}^{2}YZT+X{Y}^{2}ZT+XY{Z}^{2}T+XYZ{T}^{2} )c+ ({X}
^{2}{Y}^{2}Z+{X}^{2}Y{T}^{2}+X{Z}^{2}{T}^{2}+{Y}^{2}{Z}^{2}T )d
}$

\noindent
${ + 
({X}^{2}{Y}^{2}T+{X}^{2}Z{T}^{2}+X{Y}^{2}{Z}^{2}+Y{Z}^{2}{T}^{2}
 )e+ ({X}^{2}Y{Z}^{2}+{X}^{2}{Z}^{2}T+X{Y}^{2}{T}^{2}+{Y}^{
2}Z{T}^{2} )f
.}$

One can show that for the values of the parameters 
$$
a=u^2, \,
b=-{\frac {{u}^{2}-u+1}{2\,{u}^{2}-4\,u+1}}, \,
c= -{\frac {34\,{u}^{2}-18\,u-7}{29\,{u}^{2}+22\,u-33}}
,$$
$$
d = {\frac {7\,{u}^{2}+4\,u-6}{3\,u-2}},  \,
e= -{\frac {3\,{u}^{2}+6\,u-8}{3\,{u}^{2}+u-2}}, \,
f= -{\frac {225\,{u}^{2}-156\,u-10}{5\,{u}^{2}+212\,u-163}},$$
where $u^3+u^2-1=0 $,
the surface has the required 
singularities (see [{\bf CG}]).
Also, the $\sigma-$invariance implies that $F_5$ vanishes at $Q_0$
and is identically zero along $r$.

Let $S$ be the surface in ${\bb P}^3$ given by the equation $F_5 = 0$ and let 
$\pi: V\to S$ be a minimal resolution of singularities. Let $E_i = 
\pi^{-1}(a_i)$. By the adjunction formula,
$$K_V = \pi^*(H)-E_1-E_2-E_3-E_4,\eqno (2.1)$$
where $H$ is a hyperplane section of the quintic. 
Since the points $a_i$ are not coplanar,  
$$p_g(V) = h^0(K_V) = 0.$$ 
Now one 
checks that the pencil of quadrics 
$$\lambda YT+\mu XZ = 0 \eqno (2.2)$$
defines the 
bicanonical 
linear system on $V$. It follows from the classification of algebraic surfaces that 
$V$ is of general type.
Noether's inequality together with $p_g (V) =0$ imply that
 the irregularity $q(V)$ of $V$ is also equal to zero.

We also see that the pencil $\lambda YT+\mu XZ = 0$
has four smooth base points, plus the
singular points of the quintic.  Thus the bicanonical system
has no fixed part on the canonical model, therefore no fixed
part on the minimal model.  This implies that $V$ is minimal.  

 We have 
$$E_i ^2 = -1, \quad  E_i\cdot K_V = 1.\eqno (2.3)$$
This gives
$$K_V^2 = 5-4 = 1,$$
thus $V$ is a numerical Godeaux surface.

\vglue .5in
\centerline{{\bf 3. From quintic to double plane}} 
\bigskip
In this section we shall find a double plane model of the surface $V$.

Let 
$$R = \pi^{-1}(r) \eqno (3.1)$$
be the inverse image of the line lying on $S$ which is fixed by the involution $\sigma^2$.
This is a smooth rational curve on $V$. Since the line $r$ does not contain 
the points $a_i$, we have
$$K_V\cdot R = 1, \quad R^2 = -3.\eqno (3.2)$$

\medskip
Consider the linear system
$$ |\pi^*(H)-R|.\eqno (3.3)$$
This is the pre-image of the pencil of plane quartic curves cut out on $S$ by planes 
through the line $r$.
We have
$$\left( \pi^{\ast}H-R \right)\cdot K_V = 4, \quad \left( \pi^{\ast}H-R \right)^2 = 0.$$
Thus $| \pi^{\ast}H-R | $ is a pencil of curves of genus 3. 

The line $r'$ intersects $S$ at 5 points. Let
$$\Sigma = \{Q_0, Q_1, Q_2, Q_3, Q_4\}$$
be the set of the corresponding points on $V$. They are isolated fixed points of
$\sigma^2$ on
$V$.  The point 
$Q_0$ is fixed by $\sigma$. The other two fixed points of $\sigma$ on $V$ lie on
$R$. 

Let $a:V'\to V$ be the blowing up of the points $Q_0,Q_1,\ldots,Q_4$. The involution
$\sigma^2$ acts on 
$V'$ with the set of fixed points equal to the union of the curve 
$$R' = a^{-1}(R)\eqno (3.4)$$
and the
curves 
$$ R_i = a^{-1}(Q_i), \, i = 0,\ldots,4.\eqno (3.5)$$
Let
$$p:V'\to F = V'/(\sigma^2)\eqno (3.6)$$
be the canonical projection. The surface $F$ is nonsingular. 

\medskip\noindent
\proclaim Proposition 3.1. $F$ is a rational surface. The pencil $|\pi^{\ast}H-R|$ is the 
pre-image under $p:V'\to F$ of a pencil $|B|$ of elliptic curves
on $F$.  

{\sl Proof.} We have 
$$K_{V'} = p^*(K_F)+\sum_{i=0}^4 R_i + R'.$$
This implies that the image of $p^*(2K_F)$ on $V$ is equal to 
$2K_{V}-2R$. But the linear system $|2K_V-2R| $ is empty since no quadric from the pencil
$\lambda YT+\mu XZ = 0$ contains the line $r$. This implies that $|2K_F|=\emptyset$, hence 
$F$ is rational.

Each fiber in the pencil $|\pi^{\ast}H-R|$ is invariant with respect to $\sigma^2$. 
Its general member intersects the ramification locus of the projection $p:V'\to F$ 
at four points  lying on $R$, so
the quotient $\left(\pi^{\ast}H-R\right)/\sigma^2$ is a curve of genus 1. 
Thus $|\pi^{\ast}H-R|$ defines an elliptic pencil
$|B|$ on $F$.

\bigskip

Let us compute the topological 
Euler-Poincar\`e characteristic $e(F)$ of the surface $F$. By Noether's formula $12(1-q+p_g)
= K^2+e$, we have
$e(V) = 11.$ Thus $e(V') = 16$ and the formula
$e(V') = 2e(F)-e(R_0+\ldots+R_4+R')$ gives
$$e(F) = 14.\eqno (3.7)$$
Since $p_g(F) = q(F) = 0$, we must have $K_F^2 = -2$.
 Now observe that the elliptic curves $E_i$ on $V$ form two orbits with respect to 
$\sigma^2$. Let $\bar E_1$ and $\bar E_2$ be the two elliptic curves on $F$ representing
these two orbits. Since each $E_i$ does not intersect a general member of 
$|\pi^{\ast}H-R|$, 
the curves $\bar E_1, \bar E_2$ are 
contained in fibres of $|B|$. Also ${\bar E_i}^2 = -1$. Since for any 
minimal elliptic surface $X$ we have $K_X^2 = 0$, a relative minimal model of $F$ is an
elliptic  surface $F'$ with $e(F') = 12$. This implies that
$F$ is obtained from $F'$
by blowing up two points. Let $B_1'$ and $B_2'$ be the exceptional curves. The corresponding
two elliptic fibres of
$F$ are the divisors 
$B_1 = \bar E_1+B_1', B_2 = \bar E_2+B_2'$.  They correspond to the plane sections $X+Z = 0$ and
$Y+T = 0$  of the quintic surface. 
\bigskip
We next consider the linear system
$|3K_V-R|.$

\medskip
\proclaim Proposition 3.2. 
$|3K_V-R|$ is a pencil of curves of genus 2, without fixed part.

{\sl Proof.}
By Riemann-Roch, dim $H^0(V, {\cal O}_V(3 K_V) ) =4$. The tricanonical linear system defines
a $\sigma^2$-equivariant rational map from $V$ to ${\bb P}^3$. Since $\sigma^2$ is
identical on $R$, the image of $R$ is contained in the fixed locus of $\sigma^2$ in
${\bb P}^3$. This fixed locus is the union of two skew lines or the union of a plane and a
point.   Since the preimage of this fixed locus on $V$ is equal to $R$ plus five isolated
points, we see that the fixed locus in ${\bf P}^3$ cannot be a plane.  Also $R$ is not
mapped to a point, thus
 the image of $R$ is a line and hence the restriction of
$|3K_V|$ to $R$ is a pencil. By considering the exact sequence
$$0\to {\cal O}_V(3K_V-R)\to {\cal O}_V(3K_V)\to {\cal O}_R(3K_V)\to 0$$
we obtain that $h^0(3K_V-R) = 2$.
Thus $|3K_V - R|$ is a pencil.

It follows from (3.2) that 
$$(3K_V-R)^2 = 9-6-3 = 0,  \;   (3K_V-R)\cdot K_V = 3-1 = 2;$$
by the adjunction formula,  the genus of the pencil is two.

We next show that this pencil has no fixed part.  
Assume there is a fixed part $F$, with 
$$\left| 3K_V -R \right| = \left| M \right| +F.  \eqno(3.8) $$
Intersecting both sides of (3.8)
with $K_V$, we obtain $K_V\cdot M \le 2$. The moving part $|M|$ must have genus two,
so
$M \cdot K_V =1$ and $M^2=1$, or
$M \cdot K_V =2$ and $M^2 = 0$.  Intersecting both sides of (3.8)
with $M$ gives
$$3K_V\cdot M - M \cdot R = M^2 + M \cdot F.$$
Since $|M|$ is obviously $\sigma^2$-invariant, and $\sigma^2$ has isolated fixed
points, each member of $|M|$ is $\sigma^2$-invariant. Since its genus is $2$,
$\sigma^2$ has 2 or 6 fixed points on a general member. This gives $M\cdot R = 6$ or
$M\cdot R = 2$. The first case immediately gives $M^2 = 0$ and 
$M\cdot F = F^2 = 0$.  
This contradicts the assumption that $F\ne 0$. Assume $M\cdot R = 2$. Since $M\cdot F >
0$, this is possible only if 
$$M^2 = 0,\quad  M\cdot F =
4,\quad F^2 = -8,\quad F\cdot K_V = 0,\quad F\cdot R = 4.$$
It is easy to see that this implies that $F$ is the sum of four disjoint $(-2)$-curves
each intersecting $M$ with multiplicity $1$. None of them is
$\sigma^2$-invariant (otherwise the intersection point with each member of $|M|$ is
$\sigma^2$-invariant). Thus $\sigma$ acts transitively on the components of $F$, so
that none of them contains a fixed point of $\sigma^2$. However this contradicts the
fact that
$F\cdot R > 0$.

\medskip\noindent
\proclaim Lemma 3.3. Let $\tau$ be an involution on a smooth surface $X$,
 and let $W$ be a
one-dimensional irreducible component of the locus of fixed points of $\tau$. Let $C \ne W$
be an irreducible $\tau$-invariant curve on $X$ which intersects $W$ at a point $x$. Then
either
${\rm mult}_x(C,W) = 1$ or $x$ is a singular point of $C$. If $x$ is an ordinary node
of $C$, then ${\rm mult}_x(C,W) = 2$. If $x$ is an ordinary cusp point of $C$ then
${\rm mult}_x(C,W) = 3$.

{\sl Proof.} Since $\tau$ acts identically on $W$, it acts identically on the tangent
space $T(W)_x$. If $x$ is a nonsingular point of $C$ and $C$ is tangent to $W$ at $x$, then
$\tau$ acts on the tangent space $T(C)_x = T(W)_x$ identically. This of course implies that
$\tau$ acts identically on $C$. However this contradicts the fact that the locus of fixed
points of $\tau$ is smooth. Assume that $C$ has a node; then $\tau$ switches the tangent
directions at the two branches, so each branch intersects $W$ with multiplicity $1$. If $C$
has a cusp at $x$, then the cuspidal tangent direction coincides with the tangent direction
of $W$ at $x$, so ${\rm mult}_x(C,W) = 3$.

\proclaim Lemma 3.4. Let $\sigma$ be an involution on a nonsingular surface $X$. Let
$C$ be a $\sigma$-invariant curve on $X$ on which $\sigma$ does not act
identically. Assume that $C$ passes through $k$ isolated fixed
points $p_i$ of
$\sigma$ on
$X$. Then
$$C^2 \equiv \sum\limits_{i=1}^k({\rm mult}_{p_i}C)^2 \quad {\rm mod }\quad 2.$$

{\sl Proof.} Let $X'\to X$ be the blowing up of isolated fixed points of $\sigma$ and $\bar
C$ be its proper inverse transform.  Let $\bar X =
X'/(\sigma)$. This is a nonsingular surface and $\bar C$ is equal to the pre-image of a
curve on it not contained in the branched locus of the projection $X'\to \bar X$.
Therefore
$\displaystyle{ \bar C^2 = C^2-\Sigma_{i=1}^k({\rm mult}_{p_i}C)^2}$ is even.

\medskip\noindent
\proclaim Lemma 3.5. Let $\Sigma = \{Q_0,\ldots,Q_4\}$ be the set of isolated fixed
points of $\sigma^2$. Let $D$ be a member of the pencil
$|3K_V-R|$ such that $\Sigma\cap D\ne \emptyset$.
Then one of the following cases occurs:
\item{(i)}$\#\Sigma\cap D = 1, D = A_1+A_2$
is the union of two irreducible curves of arithmetic genus 1 with $A_1\cap A_2 = \{Q_i\}$;
\item{(ii)}$\#\Sigma\cap D = 1, D = 2A+N_1+N_2$, where
$A$ is an irreducible curve of arithmetic genus 1 passing through $Q_i\in \Sigma$ and $N_1 =
\sigma^2(N_2)$ are two disjoint $(-2)$-curves intersecting $A$ with multiplicity $1$;
\item{(iii)}$\#\Sigma\cap D = 2, D=A_1+A_2+N$ where $A_1,A_2$ are two disjoint
irreducible curves
of arithmetic genus $1$, each passing through one of the two points of $\Sigma \cap D$, and
$N$ is a $(-2)$-curve intersecting $A_1$ and $A_2$ each with multiplicity one
at these points.

{\sl Proof.}  Since $D\cdot K_V = (3K_V-R)\cdot K_V = 2$, we can write
$D = A+N$, where
$N$ is the union of $(-2)$-curves and $K_V\cdot A = 2$. The latter
equality shows that $A$ has at most two irreducible components. For
any component
$C$ of
$N$, we have $C\cdot (3K_V-R) = -C\cdot R$, therefore $C \cdot R =0$ and
 $C \cdot (3K_V -R) =0$.
Thus
$$R\cdot A = R\cdot D = R\cdot (3K_V-R) = 6.$$

If $A$ is an irreducible curve of arithmetic genus  2,
 then $A^2 = 0$ and $D = A$.  
By Lemma 3.3, $D$ has either six fixed points of
$\sigma^2$ on $R$, or $D$ has a node along $R$
and four additional fixed points on $R$, or two nodes and two additional fixed points, or
a cusp along $R$ and three additional
fixed points from $R$, or 2 cusps and no fixed points on $R$. Since $D$ has  at
least one fixed point away from $R$ it is easy to see that $\sigma^2$ acts on the
normalization of $D$ with an odd number of fixed points. This obviously contradicts the
Hurwitz formula. 

A similar argument shows that $D$ cannot be irreducible and reduced of arithmetic genus
less than $2$.

Thus $A$ is reducible or multiple.
Observe that
$D$, and hence
$A$ and
$N$, are
$\sigma^2$-invariant.  Since
a nonsingular rational component of $A$ has at most two fixed points, the equality
$A\cdot R = 6$ shows that $A$ cannot consist of two nonsingular rational components.
Also if $C$ is a nonsingular rational component of $A$, then the remaining
component $C'$ must satisfy $K_V\cdot C' = 1$ and hence must be of arithmetic genus $1$.
Since $R\cdot (C+C') = 6$ the rational component $C$ must have two distinct intersection
points with $R$ (Lemma 3.3) and no more fixed points of $\sigma$. Since $C^2 = -3$ is odd
this contradicts Lemma 3.4. 

Thus each
irreducible component of $A$ is of arithmetic genus 1. Notice that since $A\cdot R = 6$,
we must have two irreducible components $A_1$ and $A_2$ which of course may coincide.
Obviously, each $A_i$ is $\sigma^2$-invariant. By Lemma 3.4 each $A_i$ contains an odd number
of nonsingular points from $\Sigma$. Since $(A_1+A_2)\cdot R = 6$ and $A_i\cdot R\le 4$, we
obtain that each $A_i$ contains exactly one nonsingular point from the set $\Sigma$.

\smallskip
Case 1: $A_1\ne A_2, A_1\cdot A_2 \ne 0$. In this case $A_1$ and $A_2$ have a common point
from $\Sigma$ and $A_1\cdot A_2 = 1$. The formulae
$0 = D\cdot A_i = (A_1+A_2+N)\cdot A_i = -1+1+N\cdot A_i$
imply that $N = 0$. This is  case (i).

\smallskip
Case 2: $A_1\ne A_2, A_1\cdot A_2 = 0$. As above we obtain that each $A_i$ intersects
$N$ at one point. This point must be one of the points $Q_i$. Since $0 =
(A_1+A_2+D_2)^2 = 2+N^2$, we see that $N^2 = -2$.
Let $N_1,\ldots,N_k$ be the irreducible components of $N$ such that $N_1$ intersects
$A_1$ and $N_k$ intersects $A_2$. Clearly, $N_1$ and $N_k$ are 
$\sigma^2$-invariant.  If
$k > 1$,i.e. $N_1 \ne N_k$, each $N_i$ contains an additional point from
$\Sigma$.  We know that the image of $\Sigma$ in $S$ lies on the line $r'$. 
Since $N \cdot K_V =0$, $N_1$ and $N_k$ must be components of a member of the pencil
$|2K_V|$; this system corresponds to quadrics in ${\bb P}^3$, which can intersect the
line $r'$ in at most two points.  Therefore $N_1\cup N_k$ must have at most two points
from $\Sigma$. This contradiction shows that $N_1 = N_k$ and hence $N$ is a single
$(-2)$-curve. This gives us case (iii). 

\smallskip
Case 3: $A_1 = A_2$ so that $3K_V-R = 2A_1+N$. Intersecting both sides with $K_V$ we get 
$A_1\cdot K_V = 1, A_1^2 = -1$.
 The formulae
$$0 = D\cdot A_1 = (2A_1+N)\cdot A_1 = -2+N\cdot A_1,$$
$$D^2 = -4+4A_1\cdot N+N^2 = 0$$
show that $N\cdot A_1 = 2, N^2 = -4$. The last equality shows that $N$ consists of at
least two irreducible components. It cannot consist of two components $N_1$ and $N_2$
with $N_1\cdot N_2 = 1$. Indeed, write $N = mN_1+nN_2$ to get 
$$N^2 = -2m^2-2n^2+2mn = -2(m-n)^2-2mn.$$
It is easy to see that $N^2$ cannot be equal to $-4$. So, $N = N_1+N_2$, where $N_1\cdot
N_2 = 0$ as in case (ii), or $N$ contains more than two irreducible components. It remains
to exclude this case. Observe that $D$ is not
$\sigma$-invariant since otherwise
$A_1$ is
$\sigma$-invariant and hence has two fixed points of $\sigma$ (it is an easy fact that
which follows from the Hurwitz formula that any involution of order 4 on an irreducible
curve of arithmetic genus 1 has 2 fixed points). Let $D_1,\ldots,D_k$ be the divisors from
$|3K_V-R|$ which contain points of
$\Sigma$ and let 
$r_i$ be the number of irreducible components of $D_i$. Since
$b_2(V) = {\rm rk\  Pic}(V) = 9$, we have $\sum_i(r_i-1)\le 7$. If our $N$ contains more
than three components, $D$ and $\sigma(D)$ will contribute at least $4+4$ to this sum which is
too much.  Since
$D$ is
$\sigma^2$-invariant, $N$ must be supported on the nodal cycle $N_1+N_2+N_3$ of type
$A_3$ with $N_1$ and $N_3$ intersecting $A$. Moreover $N_2$ must be
$\sigma^2$-invariant and hence contains two points from
$\Sigma$. Thus $D$ contains three points from $\Sigma$, and $\sigma(D)$ contains three more
points from $\Sigma$ which is too many.

This proves the assertion.

\medskip \noindent
\proclaim  Corollary 3.6. In the notation of the previous lemma, let $D_i, i =
1,\ldots,k,$ be the set of divisors from $\left| 3K_V-R \right|$ 
such that $D_i\cap \Sigma \ne
\emptyset.$ Then one of the following cases occurs:
\item{(i)} $k = 5$, each $D_i$ is of type (i);
\item{(ii)} $k = 5$, $D_4$ and $D_5 = \sigma(D_4)$ are of type (ii), $D_1,D_2,D_3$ are of
type (i); 
\item{(iii)} $k = 3$, $D_2$ and $D_3 =\sigma(D_1)$ are of type (iii), $D_1$ is of type
(i).

\noindent
{\sl Moreover, the unique $\sigma$-invariant divisor, say $D_1$, is of type (i)}. 

{\sl Proof. } Since
$\Sigma$ consists of five points, and each divisor $D_i$ contains at most two points from $\Sigma$, we have $5\ge k\ge 3.$
Assume
$k = 5$. Then each
$D_i$ is of type (i) or (ii). Let
$r_i$ be the number of irreducible components of $D_i$. As we already noticed in the
previous proof,  $\sum_i(r_i-1)\le 7$. This implies that we
have at most two divisors $D_i$ of type (ii). If there is only one $D_i$ of type (ii),
then it is necessary $\sigma$-invariant. But then each of the rational components of
$D_i$ is $\sigma^2$-invariant. Since it does not intersect $R$ and does not contain
points from $\Sigma$, we get a contradiction. This gives us two cases (i) and (ii)
with $k = 5$.

Assume $k = 4$. Then exactly one of the divisors $D_i$ must contain two points from
$\Sigma$, and hence it must be of type (iii) and $\sigma$-invariant. But then its
rational component must be invariant and contain two isolated fixed points of
$\sigma$. This is impossible since we have only one such point. So this case does not
occur.

Assume $k = 3$. Then we must have two divisors of type (iii). The third divisor cannot
be of type (ii) since otherwise it is $\sigma$-invariant, and we have already seen that
this is impossible.

\medskip \noindent 
\proclaim Theorem 3.7. $V$ is isomorphic to a nonsingular minimal model of the double
cover of ${\bb P}^2$.   The  branch curve is an
irreducible curve
$W$ of degree 10 with a point
$q$ of multiplicity 4 and 10 triple  
points $p_i, p_i'$, $i = 0,\ldots,4$, where $p_i'$ is infinitely near to $p_i$. The points 
$p_0,\ldots,p_4,q$  do not lie on a conic.  
Also there exists a pencil $|\bar B|$ of elliptic curves of degree 4 with double points
at 
$p_0$ and 
$p_0'$ and simple points at other triple points of $W$.

{\sl Proof. } We shall employ the notation from the proof of Proposition 3.1. Let $B$
be the pencil of elliptic curves on $F = V'/(\sigma^2)$ whose general fibre is the
quotient of the general fibre of $|\pi^*(H)-R|$ by $(\sigma ^2)$. Let
$|D|$ be the pencil of rational curves on $F$ whose general fibre is the
quotient of the general fibre of $|3K_V-R|$ by $(\sigma ^2)$. It folows from (3.2)
that $(3K_V-R)\cdot (\pi^*(H)-R) = 8$, hence 
$$B\cdot D = {1\over 2}(3K_V-R)\cdot (\pi^*(H)-R) = 4.$$
Let $B_1 = \bar E_1+B_1'\in |B|$ be one of the two non-minimal fibres of the elliptic
fibration defined by $|B|$. The curve $B_1'$ is the exceptional curve of the first
kind and the pre-image of $\bar E_1$ is the union of two elliptic curves $E_1$ and
$E_3$ on $V$. Since $(E_1+E_3)\cdot (3K_V-R) = 6$, we get 
$$ D \cdot \bar E_1 = 3, \quad  D\cdot \bar B_1' = 1.$$
Thus $B_1'$ is a section of the pencil $|D|$. 

Now consider the divisors $D_i, i= 1,\ldots,k,$ from $|3K_V-R|$ which contain a point
from the set $\Sigma$ of isolated fixed points 
of $\sigma^2$. Lemma 3.5 and  Corollary 3.6 describe these divisors.

Let $\bar D_i$ denote the corresponding divisor in $|D|$.  If $D_i = A_1+A_2$ is
of type (i), then 
$$\bar D_i = \bar A_1^{(i)}+\bar A_2^{(i)}+\bar R^{(i)},$$
where $\bar A_1^{(i)}, \bar A_2^{(i)}$ are exceptional curves of
the first kind, and $\bar R^{(i)}$ is the image in $F$ of $a^{-1}(Q)\in V'$, where 
$Q = D_i\cap \Sigma$. 

If $D_i = 2A +N_1+N_2$ is of type (ii), then similarly we have
$$\bar D_i = 2 \bar A + \bar N^{(i)} + \bar R^{(i)}.$$
where $\bar R^{(i)}$ is the image of $a^{-1}(D_i\cap \Sigma)$ and $\bar N^{(i)}$ is the
image of the orbit of $(-2)$-curves $N_1 + N_2$ on $F$.  

Finally if $D_i$ is of type (iii),
$$\bar D_i = \bar A_1^{(i)} + \bar A_2^{(i)} + \bar R_1^{(i)} + \bar R_2^{(i)}+ \bar
N^{(i)},$$ 
with the similar notation.

\bigskip
\epsfxsize=4.5in
\epsfbox{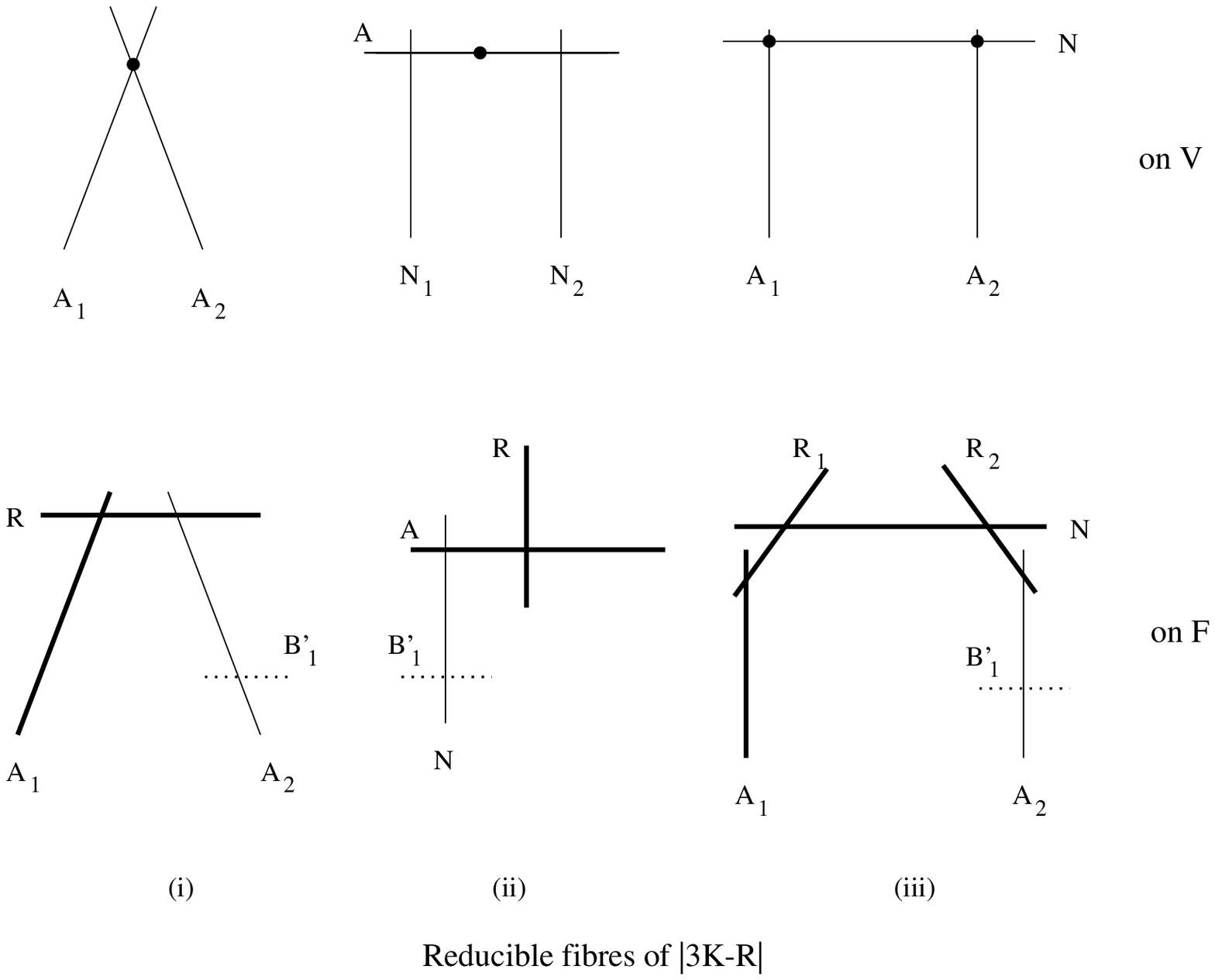}
\bigskip
Consider the section $\bar B_1'$. It intersects each divisor $\bar D_i$ simply at one
point. Without loss of generality we may assume that this point
does not lie on $\bar A_1^{(i)}$. Also we observe that $\bar B_1'$ does not intersect the
the curve $\bar R^{(i)}$ (where $\bar R^{i} = \bar R_1^{i}+\bar R_2^{(i)}$ in case $D_i$
is of type (iii)). This is because the planes
$X+Z= 0$ and $Y+T = 0$ do not contain points from $\Sigma$. 

For each $i= 1,\ldots,k$, let ${\cal E}_i$ be the divisor chosen as follows. If $D_i$ is
of type (i) or (ii),
$${\cal E}_i = \bar A_1^{(i)}+\bar R^{(i)},$$
if $D_i$ is of type (iii),
$${\cal E}_i = \bar A_1^{(i)}+\bar R_1^{(i)} + \bar R_2^{(i)}+ \bar
N^{(i)},$$
It is easy to see that starting from the curve $\bar A_1^{(i)}$ we can blow down the
divisor ${\cal E}_i$ to a nonsingular point. Using Corollary 3.6, it is easy to check
that  the union of the divisors
${\cal E}_i, i = 1,\ldots,k,$ is an exceptional divisor with 10 irreducible components
which can be blown down to $k$ distinct nonsingular points. Let
$\alpha:F\to
\bar F$ be its blowing down. It follows from (3.7) that $e(\bar F) = 4$, hence $\bar F$
is a minimal ruled surface.  Since each
${\cal E}_i$ is a part of a member of the pencil $|D|$, the
image of
$|D|$ on
$\bar F$ defines a ruling 
$$\alpha_1:\bar F\to {\bb P}^1.$$ 
It is clear that the curve  $B_1'$ does not
intersect
${\cal E}_i$ if $D_i$ is of type (i) or (ii). We shall prove a little later that case
(iii) does not occur. Let us continue the proof
accepting this fact.  The image of
$B_1'$ is a  curve
$L$ on
$\bar F$ with
$L^2 = B_1'^2 = -1$. This implies that $\bar F$ is the minimal ruled surface ${\bf F}_1.$ 
We can blow down the curve $L$ to arrive at ${\bb P}^2$.

Since we have blown down all the curves $\bar R_j^{(i)}$, our surface $V$ is birationally
isomorphic to the double plane branched along the curve $W$ equal to the image of 
the ``line'' $R$. The image $\bar R$ of $R$ on $F$ has self-intersection $-6$. It
intersects each  ${\cal E}_i$ at the curve $\bar A_1^{(i)}$ with multiplicity $3$ and
intersects the curve $L$ with multiplicity $4$. When we blow down ${\cal E}_i$ we
increase the self-intersection by $18$. When we blow down $L$ it increases by $16$.
Using Corollary~3.6 we see that in any case we increase $\bar R^2$ by 106, and hence
$W^2 = 100$, i.e. $W$ is a curve of degree 10. It has a quadruple point at the image
of $L$.  It has a triple point with one infinitely near triple point at the image of
${\cal E}_i$.  Applying Corollary 3.6 we get
the assertion about the degree and the singularities of the branch curve $W$.

It remains to check the assertion about the  elliptic pencil. As the reader can
already guess this must be the image of the pencil $|B|$. First of all, we have already
noticed that
$B\cdot D = 4$. This means that the image of
$|B|$ on
${\bf F}_1$ intersects the ruling with multiplicity $4$. This implies that the image
$\bar B$ of
$B$ on ${\bb P}^2$ is a pencil of elliptic quartics. Let $D_1$ be the unique
$\sigma$-invariant divisor among the divisors $D_1,\ldots,D_k$. We know from   
Corollary~3.6 that it is of type (i).   
A general member $B$ of the pencil $|\pi^*(H)-R|$
intersects  each $D_1 = A_1+A_2\in |3K_V-R|$ with
multiplicity 4. Since $\sigma(A_1) = A_2$ and 
$\sigma (\pi^*(H)-R) = \pi^*(H)-R$ this implies that $(\pi^*(H)-R)\cdot A_1 =
(\pi^*(H)-R)\cdot A_2 = 4$, and hence
$$(B_1'+\bar E_1)\cdot \bar A_1 = (B_1'+\bar E_1)\cdot \bar A_2 = \bar E_1\cdot \bar A_1 = 2.$$
 Thus when we blow down ${\cal
E}_1$ the image of $B$ in ${\bb P}^2$ acquires a
double point at $p_0$ and a double point at
$p_0'$. There are no more singular points for a general member of $|\bar B|$. So $B$
intersects each ${\cal E}_i$ with multiplicity 1. This gives us the assertion about
the pencil $|\bar B|$.  The assertion about the conic follows easily from computing the
canonical class of the double plane (see the proof of Theorem~4.1).   

It remains to pay our debt by checking that case (iii) in Corollary~3.6  does not
occur.  Suppose $D_i$ is of type (iii); we claim that $B_1'$ does not intersect
${\cal E}_i$ in this case. Suppose it does. Then we still can blow down the
divisors ${\cal E}_1,{\cal E}_2, {\cal E}_3$ to arrive at a minimal ruled surface $\bar
F$.  The image $\bar B_1'$ of $B_1'$ on $\bar F$ is a section with positive
self-intersection. Its divisor class on $\bar F$ is equal to  $nf+s, n \ne 0$, where $f$
is the class of the fibre and
$s$ is the class of the exceptional section. Recall that $R' = R/(\sigma^2)$ is a
multi-section of the pencil $|D|$ of degree 6. So the image $\bar R$ of $R'$ 
on $\bar F$ intersects $f$ with multiplicity $6$. This gives 
$$\bar B_1'\cdot \bar R = (nf+s)\cdot \bar R \ge 6n \ge 6.$$
However $B_1'$ is the image on $F$ of a component of a divisor from
$|\pi^*(H)-R|$ which intersects $R$ with multiplicity 4. This shows that $\bar
B_1'\cdot \bar R
\le 4$ and gives us a contradiction. So, even in the case when $D_i$ is of type (iii) we
can blow down the exceptional divisors ${\cal E}_i$ and $B_1'$ to obtain the plane
model of $F$. The image of the pencil $|B|$ is again a pencil of quartics. Thus the
self-intersection of $B$ must decrease by 16. When we blow-down ${\cal E}_1$ and $B_1'$
we decrease it by $8$. Thus ${\cal E}_2$ and ${\cal E}_3$ must contribute the
remaining 8. Since the general member of $\bar B$ has only two singular points which give
a tacnode at the image of ${\cal E}_1$, we see that $B\cdot {\cal E}_2 = B\cdot {\cal
E}_3 = 1$. We also know that $B$ does not intersect the components $R_j^{(i)}$. This
easily gives that $B^2$ decreases by at most $3+3 = 6$, a contradiction.   Thus, all
the assertions of the theorem have been checked, and we obtain that 
case (iii) in Corollary~3.6 does not occur.  

\medskip\noindent
{\bf Remark 3.8.}   
Recall that  the elliptic pencil $|B|$ on $F$ contains two reducible members $\bar
E_1+B_1'$ and $\bar E_2+B_2'$. The image of 
$\bar E_1$ is a quartic from the pencil
$|\bar B|$ passing through
$q$. The curve $B_1'$ is of course blown down to $q$.
The image of $\bar E_2$, the
orbit of $E_2$ and $E_4$,
 is a cubic curve which passes simply through the points 
$p_i,p_i', i = 0,\ldots,4$. The image of $B_2'$ is the line passing through $p_0$ and
$p_0'$.   One can make a standard (degenerate) Cremona transformation with center at 
the points $p_0,p_0',q$. It transforms the branch curve $R$ into itself, 
and replaces the curve  $L$ with the line $l$ through $(p_0,p_0')$. This is induced by
the involution 
$\bar \sigma$.

\medskip\noindent
{\bf Remark 3.9.} We shall later show that also the case (ii) in Corollary 3.7 does not
occur for the surfaces which we are considering.  It is quite possible that  cases
(ii) and (iii) occur for some genus 2 pencil on a numerical Godeaux surface.
Observe that  the case (iii) corresponds to the degeneration when, say
$p_1 = p_2$ and
$p_3 = p_4$. The case (ii) corresponds to the degeneration when the points
$q,p_i,p_i'$ become collinear.  
.

\medskip\noindent
{\bf Remark 3.10.}  
We can also  realize $V$ as a double cover of ${\bb P}^1 \times {\bb P}^1$ as follows.  
Consider the linear system $|4K_V -R|$;  we have
$$(4K_V-R)^2 = 5, \, (4K_V-R) \cdot K_V = 3, \, 
(4K_V-R) \cdot (3K_V-R) =2.$$
Thus $|4K_V-R|$ is a linear system of genus five curves.   
One checks that this  is a pencil with five base points, which are exactly the
points
$Q_0, Q_1,
\dots, Q_4$ fixed by the involution $\sigma^2$.
We know that the members of $|3K_V-R|$ passing through these points contain pairs
of elliptic curves $A_i +A_i'$.

Thus these two pencils define a double cover $V \rightarrow {\bb P}^1 \times {\bb P}^1$.
This map blows up the base points of $|4K_V-R|$ to rulings of ${\bb P}^1 \times {\bb P}^1$
and contracts the elliptic curves to pairs of points on these fibres (for details, see
[{\bf R}]).  
In case (ii) of the
Lemma,  we have $A_i =  A_i'$ for a member of $|3K_V -R|$;
in this case the pair
of  points on the corresponding fibre are infinitely near.

We have 
$$(3K_V-R) \cdot R = 6, \; (4K_V-R) \cdot R =7,$$
thus the rational curve $R$ is mapped to a curve of bidegree $(6,7)$ with pairs of
triple points on the five special fibres of $|3K_V -R|$. This curve on ${\bb
P}^1 \times {\bb P}^1$ can then be birationally transformed into 
a degree ten plane curve as in the theorem.

\vglue .5in
\centerline{{\bf 4.  From double plane to quintic}}
\medskip
 One can reverse the previous
construction.

\medskip 
\proclaim Theorem 4.1. Assume there exists an irreducible plane curve $W$ of degree 
10 with an point $q$ of multiplicity 4 and 10 triple 
singular 
points $p_i,p_i', i= 0,\ldots, 4$, where $p_i'$ is infinitely near to $p_i$. We assume
that  the points 
$p_0,\ldots,p_4,q$ do not lie on a conic.  
Also assume that there exists a pencil $D$ of elliptic curves of degree 4 with double points at $p_0$ and 
$p_0'$ and simple points at $p_i,p_i', i\ne 0$. Let $V'$ be a double cover of ${\bb P}^2$
branched over 
$W$. Then its nonsingular minimal model $V$ is a numerical Godeaux surface. Let $E$ be the
pre-image of  the elliptic pencil on $V$ and $R$ be the proper inverse transform of the curve
$W$; then the linear system
$|E+R|$ maps $V$ onto a normal quintic in ${\bb P}^3$ with 4 simple 
elliptic 
singular points of degree 1. The singular points of the quintic are the images of the inverse 
transforms of the quartic $Q$ from the pencil which passes through $q$ and the 
cubic curve $C$ which together with the line 
$l =<p_0,p_0'>$ forms a member of the 
pencil of elliptic curves.
If additionally, $W$ is invariant with respect to the standard Cremona transformation with
center 
 $p_0,p_0',q$, then the quintic is invariant with respect to a projective transformation
of order 4.

{\sl Proof.} It is rather standard. Let $\pi:F\to {\bb P}^2$ be the blowing-up of points 
$q,p_i,p_i',$ $i= 0,\ldots, 4$. Then 
$$K_F = -3h+Z+\sum_{i=0}^4(Z_i+2Z_i'),\eqno (4.1)$$
where $h$ is the pre-image of a line, $Z_i+Z_i' = \pi^{-1}(p_i), Z_i'=\pi^{-1}(p_i'), Z=\pi^{-1}(q)$.
The proper transform of $W$ is
$$\bar W = 10h-4Z-\sum_{i=0}^4(3Z_i+6Z_i').\eqno (4.2)$$
We take the double cover $p:V'\to F$ branched along the nonsingular divisor
$$W' = \bar W+\sum_{i=0}^4 Z_i.\eqno (4.3)$$
Then, we have
$$K_{V'} = p^*(K_F+{1\over 2}W') =2 p^*(h)-p^*(Z)-\sum_{i=0}^4 p^*(Z_i').$$
It follows that
$$p_*({\cal O}_{V'}(K_{V'})) = {\cal O}_{F}(K_F+{1\over 2}W)\oplus 
{\cal O}_{F}(K_F) ={\cal O}_{F}(2h-Z-\sum_{i=0}^4 Z_i' )\oplus {\cal O}_{F}(K_F).$$
Since $h^0(2h-Z-\sum_{i=0}^4 Z_i') = 0$ (otherwise there would be a conic through the points 
$q,p_i$), we obtain 
$$h^0(K_{V'}) = h^0(p_*({\cal O}_{V'}(K_{V'}))) = 0.$$
By Riemann-Roch, 
$h^1(2h-Z-\sum_{i=0}^4 Z_i') = 0$, so we obtain
$$h^1(K_{V'}) = h^1(p_*({\cal O}_{V'}(K_{V'}))) = 0.$$

On the other hand, 
$$|2K_{V'}| = |p^*(4h-2Z-2\sum_{i=0}^4 Z_i')|\ne \emptyset,$$
since it contains the positive divisors equal to the sum of the curves $Z_i$ and the proper 
transform of a quartic curve $C_4$ from the pencil of quartics passing simply through the 
points $p_i,p_i'$ and with multiplicity 2 at $q$. 
The curves $\bar A_i = p^*(Z_i')$ are elliptic curves with self-intersection $-2$.
The curves $R_i = p^*(Z_i)$ are exceptional curves of the first kind. Blowing down the curves 
$R_i$ we get a surface $V$ with 
$$K_V^2 = {K_{V'}}^2+5 = 1 .$$ 
Together with the property $h^0(2K_V) \ne 0$, this shows that $V$ is a surface of general type with 
$$p_g = q = 0.$$
It is also minimal. In fact any exceptional curve of the first kind intersects
$2K_V$ negatively, 
hence must be a component of 
the image on $V$ of the proper transform of any quartic curve $C_4$ from above. It is easy to see that the pencil of such quartics does not have base components.
Thus $V$ is a numerical Godeaux surface. 

Let $C_1 \sim h-Z_0-2Z_0'$ be the proper transform on $F$ of the line $l$, and 
$C_2 = Z$. 
The pre-image on $F$ of the 
elliptic pencil formed by the quartic curves is equal to
$$4h-2Z_0-4Z_0'-\sum_{i=1}^4(Z_i+2Z_i');$$
we have 
$$4h-2Z_0-4Z_0'-\sum_{i=1}^4(Z_i+2Z_i')\;  \sim \;
C_1+[ 3h-\sum_{i=0}^4 (Z_i+2Z_i') ]$$
and
$$4h-2Z_0-4Z_0'-\sum_{i=1}^4(Z_i+2Z_i')\;  \sim \;
 C_2+[ 4h-2Z_0-4Z_0'-\sum_{i=1}^4(Z_i+2Z_i')-Z ].$$
The curves $B_1 \in |3h-\sum_{i=0}^4(Z_i+2Z_i')|$ and 
$B_2\in |4h-2Z_0-4Z_0'-\sum_{i=1}^4(Z_i+2Z_i')-Z|$ are elliptic curve on $F$ of self-intersection 
$-1$ which are disjoint from the branch divisor $W'$. They are equal to the 
proper transforms of the cubic curve $C$ and the quartic $Q$, respectively.
Under the double cover $p:V'\to F$ they are split into the disjoint 
sum of two elliptic curves with self-intersection $-1 \,$,
$p^*(B_1) = E_1'+E_2'$ and $ p^*(B_2) = E_3'+E_4'.$ 
We have  
$$K_{V'}+E_1'+E_2'+E_3'+E_4' \sim
[p^{\ast} ( 2h - Z -\sum_0^4 Z_i')] +$$
$$+[p^{\ast} (4h -Z -2Z_0 -4Z_0' - \sum_1^4 (Z_i +2Z_i'))]
+[p^{\ast}( 3h  -\sum_0^4 (Z_i +2Z_i'))]$$
$$=9p^*(h)-2p^*(Z)-3p^*(Z_0)-7p^*(Z_0')-2\sum_{i=1}^4 p^*(Z_i)-5\sum_{i=1}^4 p^*(Z_i').$$
On the other hand, we have
$$ p^*(4h-2Z_0-4Z_0'-\sum_{i=1}^4(Z_i+2Z_i'))+{1\over 2}p^*(\bar W+\sum_{i=0}^4 Z_i)=$$
$$= p^*(4h-2Z_0-4Z_0'-\sum_{i=1}^4 (Z_i+2Z_i'))+p^*(5h-2Z-\sum_{i=0}^4 (Z_i+3Z_i')).$$
Comparing, we get that 
$$K_{V'}+E_1'+E_2'+E_3'+E_4'\sim p^*(D')+{1\over 2}p^*(\bar W+\sum_{i=0}^4 Z_i),$$
where $|D'| = |p^*(4h-2Z_0-4Z_0'-\sum_{i=1}^4 (Z_i+2Z_i'))|$ is the pencil of curves 
originating from the pencil of quartic curves.  Blowing down the exceptional curves 
$p^{\ast} Z_i$, we get
$$K_{V}+E_1+E_2+E_3+E_4\sim R+D,$$
where $E_i,D$ are the images of $E_i',D'$, and $R$ is the reduced pre-image of $\bar W$.
Clearly, $R$ is a nonsingular model of our branch curve of degree 10. Since the genus of the
latter is zero, $R$ must be isomorphic to ${\bb P}^1$. Let us denote the linear system
$|R+D|$ by $|H|$ and prove that $|H|$ maps $V$ birationally  onto a normal quintic in ${\bb
P}^3$.

We immediately check that $D\cdot R = 4$, hence $|D|$ is a pencil of curves of genus 3. Since
$D\cdot D = 0$, we obtain $$D\cdot K_V = 4.$$ From this it follows that 
$$5 =
K_V^2+K_V\cdot(E_1+\ldots + E_4) = K_V\cdot D+K_V\cdot R = 4+ K_V\cdot R,$$
hence 
$$K_V\cdot
R = 1,\quad R^2 = -2-1 = -3.$$
 This immediately implies
$$\quad H^2 = (D+R)^2 = 5, \quad H\cdot K_V = 5$$ 
Since we know that 
$p_g(V) = q(V) = 0$, by Riemann-Roch
$$h^0(H) = {1\over 2}(H^2-K_V\cdot H)+1+h^1(H) = 1+h^1(H) = 1+h^1(K_V-H).$$
Now the exact sequence 
$$0\to {\cal
O}_V(-E_1-E_2-E_3-E_4)\to {\cal O}_V\to {\cal O}_{E_1+E_2+E_3+E_4} \to 0$$ 
together with the isomorphism ${\cal
O}_V(-E_1-E_2-E_3-E_4)\cong {\cal O}_V(K_V-H)$ gives 
$$h^1(H) = h^1(-E_1-E_2-E_3-E_4) = h^0({\cal
O}_{E_1+E_2+E_3+E_4})-h^0({\cal O}_V) = 3.$$ 
Thus $${\rm dim}|H| = h^0(H)-1 = 3.$$ The complete linear system
does not have base points. In fact, since $H=D+R$, where $|D|$ is a base-point-free pencil, any
base point of
$H$ belongs to $R$. But $H\cdot R = 1$ and the exact sequence 
$$0\to {\cal O}_V(D)\to {\cal O}_V(H)\to {\cal
O}_R(H)\to 0$$ 
shows easily that $|H|$ cuts out a pencil on $R$. This checks
that there are no base points on
$R$. 

Now the linear system $|H|$ defines a regular map $f:V\to {\bb P}^3$. Since $H\cdot E_i = 0$, 
it blows down the curves $E_i$. Since $H^2 = 5$ is a prime number, $f$ must be of degree 1 with the image a quintic surface with 4 elliptic singular points of degree 1.
Since  a nonsingular member of $H$ is of genus $g = 1+{1\over 2}(H^2+H\cdot K_V) = 6$, and this coincides with the arithmetic genus of a plane section of a quintic surface, the 
image of $f$ is a normal surface. 
So, all the assertions of the theorem are proven.
\vfill  
\eject
%\vglue .5in
\centerline{{\bf 5.  The fundamental group}}
\medskip
We shall start with the following

\medskip

\proclaim Lemma 5.1. Assume $V$ contains a  $(-2)$-curve $C$. Then $C$ is a 
component of one of the members of the pencil $|3K_V-R|$ which passes through a point $Q_i,
i\ne 0$. Moreover, 
its image in the quintic model is  a conic which passes through two 
of the four singular points $a_i$.

{\sl Proof.} Suppose $C$ is a $(-2)$-curve.  
Since $(3K_V-R)\cdot C = -R\cdot C$ and $C\ne R$, we get that 
$C\cdot R = 0$ and $C$ is contained in a member of $|3K_V-R|$. 
If $C$ does not contain one of the points $Q_i$, then
the involution $\sigma^2$ has no fixed points on $C$,
and the orbit of $C$ under the group $(\sigma)$ 
consists of 4 disjoint $(-2)$-curves.  
Using the argument from the proof of Corollary 3.6, we find that there
are too many irreducible components contained in reducible fibres of
$|3K_V-R|$. 

Thus $C$ must belong to one of the divisors $D_i$ from Corollary 3.6. 
As we have already noticed in the proof of Theorem 3.7, the case (iii) of Corollary
3.6 does not occur for our surfaces. Thus there are four of the curves $C$ forming an
orbit with respect to $(\sigma)$. Each $(\sigma^2)$-orbit is contained in one reducible
fibre of
$|3K_V-R|$. 

 Let
$D =2A+C+\sigma^2(C)$ be one of the reducible members of the pencil
$|3K_V-R|$ which contains a $(-2)$-curve $C$ as in
case (ii) of Lemma 3.5, and let $E_i$ be one of the four elliptic curves
blown up from the singular points of the quintic. Then
$E_i\cdot (3K_V-R) = 3$. Hence
$3 = E_i\cdot (2A+C+\sigma^2(C))= 2(E_i\cdot A)+E_i\cdot C+E_i\cdot \sigma^2(C)$. It follows
from the proof of Theorem 3.7 that $E_i\cdot A \ne 0$. This
implies that
$E_i\cdot A = 1$ and 
$E_i$ intersects $C+\sigma^2(C)$ at one point. Now repeating this argument for the member 
$\sigma(D)$ of $|3K_V-R|$ we find that $E_i$ intersects one more
$(-2)$-curve $\sigma(C)$ (or $\sigma^3(C)$). Thus each $E_i$ intersects exactly two $(-2)$-
curves. This implies that each $(-2)$-curve $C$ intersects two elliptic curves $E_i$ each with
multiplicity $1$. Now
$$0 = K_V\cdot C = (\pi^*(H)-E_1-E_2-E_3-E_4)\cdot C$$
gives $\pi^*(H)\cdot C = 2$, i.e. the image of $C$ in the quintic model is a conic passing
through two singular points $a_i$.

\medskip
\proclaim Lemma 5.2.
$$H_1(V,{\bb Z}) = 0.$$

{\sl Proof.} It is enough to show that ${\rm Tors} \, H_1(V,{\bb Z}) = 0$.
Consider the bicanonical system  $|2K_V|$ on $V$.  
A nonzero torsion divisor $\tau$ will cause a splitting 
$$2 K_V = \left( K_V + \tau \right) + \left( K_V -\tau  \right) .$$ 
For each nonzero torsion divisor $\tau$ the linear systems $|K_V + \tau|$ consists of a unique
curve $C_\tau$ with $\chi(C_\tau,{\cal O}_{C_\tau}) = -1$. 
We have $C_\tau \cdot C_{-\tau} =1$; moreover these two curves must intersect at a base point
of $|3K_V|$ ([M], Theorem~2).    
Since $|3K_V-R|$ is base point free, this intersection point must lie on $R$,  thus
this member of $|2K|$ must intersect $R$ at one point with multiplicity two. 

Recall that the bicanonical system $|2K_V|$ is generated
by the pullback of the pencil of quadrics ${\cal Q}(\lambda,\mu)= \lambda Y T +  \mu X Z$
in ${\bb P}^3$.
We can also consider this system as the pullback of a pencil on $V / \sigma^2$;
since $2 K_V \cdot R =2$, the general member of this pencil will intersect this
rational curve at two points, and there will be exactly two members which 
intersect  $R$ at one point with multiplicity two. We have already two of them 
coming from the two reducible quadrics. In fact,
 $$(\pi^*(H) - 2 E_i - E_j - E_k) \cdot R
=1,$$
and $\sigma^2$ switches the two components of the reducible quadric.
Thus any splitting of $2K_V$ as $C_{\tau} + C_{-\tau}$
must come from one of the reducible quadrics.  

Let us consider these reducible  
members $Y T$ and $X Z$
of the  bicanonical system. 
Each hyperplane section $X, Y, Z, $ and $T$ of ${\bb P}^3$ is equivalent
in $V$  to
$\pi^*(H) - 2 E_i - E_j - E_k$ for the appropriate indices
$i, j, k$.  We have
$$(\pi^*(H) - 2 E_i - E_j - E_k)^2 = -1$$ and
$$(\pi^*(H) - 2 E_i - E_j - E_k) \cdot K_V =1,$$
thus each irreducible component of the reducible quadrics defines on $V$ a curve of
arithmetic genus
$1$.   Since ${C_\tau}^2 = C_\tau \cdot K_V =1$, 
this shows that the curves $C_\tau$ cannot be irreducible.  

Suppose
$C_\tau$ is reducible.   
Since $K_V$ is nef, and
$C_\tau\cdot K_V = 1$, $C_\tau$ can contain only one irreducible
component which is
not a
$(-2)$-curve.  Let us denote this component by $C_\tau'$.
Let us write $C_\tau = C_\tau'+Z$, where $Z$ is the sum of $(-2)$-curves. Since $
C_\tau\cdot\pi^*(H) = K_V\cdot \pi^*(H) = 5$, we have $Z\cdot\pi^*(H) \le 4$. From the
previous lemma we know that for each irreducible component $C$ of $Z$, we have
$Z\cdot\pi^*(H) = 2$. Thus $Z = 0,$ or $Z = C$ or $Z = C_1+C_2$, where we do not assume
that $C_1\ne C_2$.

Case 1: $Z = C_1+C_2$. We have
$$1 = K_V\cdot E_i = C_\tau\cdot E_i = C_\tau'\cdot E_i +(C_1+C_2)\cdot E_i \ge
(C_1+C_2)\cdot E_i.\eqno(5.1)
$$
Recall that each $(-2)$-curve $C$ intersects two $E_i$'s. If $C_1 = C_2$, we take $E_i$
such that
$C_1\cdot E_i > 0$ and get a contradiction. If $C_1\ne C_2$ and $C_1,C_2$ are not in the
same divisor $D\in |3K_V-R|$, then we can find $E_i$ which intersects both $C_1$ and
$C_2$. This follows from a simple combinatorial argument using the fact 
that each
$E_i$ intersects exactly one
$(-2)$-component  of a reducible divisor $D$ from $|3K_V-R|$. So, applying (5.1), we get
a contradiction.  

If $C_1$ and $C_2$ are in the
same divisor $D\in |3K_V-R|$, then 
$$0=K_V \cdot (C_1 + C_2) = C_\tau' \cdot (C_1 + C_2) -4$$
implies that  $C_\tau' \cdot (C_1 + C_2) =4$; however this contradicts $C_{\tau}' \cdot (3K_V -R) =2$.

Case 2: $Z$ is irreducible. 
If $Z$ is an irreducible  $(-2)$-curve, then we have 
$K_V\cdot C_\tau' = 1, 
C_\tau'\cdot Z = 2, C_\tau'^2 = -1$. Thus 
$C_\tau'$ is of arithmetic genus 1 and $C_\tau'\cdot Z = 2$. In the quintic model 
$C_\tau'$ is a plane cubic. The formula (5.1)  shows
that $C_\tau'$ intersects exactly two
 $E_i$'s, say $E_1$ and $E_2$, with multiplicity 1. The residual conic defines a
divisor
$Z'\in |\pi^*(H)-C_\tau'-aE_1-bE_2|$. We have 
$$Z'^2 = -2+2a+2b-a^2-b^2, \quad Z'\cdot K_V = 4-a-b.$$ 
Since $Z'^2+Z'\cdot K_V = -2$, we obtain
$$a(a-1)+b(b-1) = 4.$$
The only solution is $a = b = 2$.  Thus
$$Z' \sim \pi^*(H)-C_\tau'-2E_1-2E_2$$
and $Z'$ is a $(-2)$-curve.
 
Since the image of $C_\tau'$ in the quintic model is a plane cubic,   
one of the plane components of the reducible quadric which contains
$\pi(C_\tau')$ must also contain the residual conic $Z'$, and 
$$\pi^*(H) - 2 E_i - E_j - E_k \sim C_\tau'+Z' \sim \pi^*(H)-2E_1-2E_2.$$
This is obviously impossible. 

Therefore a decomposition of $2K$ as $(K_V + \tau) + ( K_V - \tau)$ cannot occur,
and we have shown that there are no torsion divisors on $V$. 

\bigskip

We next prove the following result:
\medskip

\proclaim Theorem 5.3. 
$$\pi_1(V) = \{1\}.$$

\bigskip
By Lemma~5.2, it suffices to show that $\pi_1(V)$ is commutative.

As before, write $V' \to V$ for the blowup of the five isolated
fixed points of $\sigma^2$, and $p: V' \to F$ for the double cover
of the rational surface $F$ branched along the divisor
$W' = \bar W + \sum_0^4 Z_i$ (see (4.3)). Recall that $R' = p^{-1}(\bar W)$ is the 
pre-image of the curve $R$ on $V'$, and  $R_i = p^{-1}(Z_i)$ are the five $(-1)-$curves on
$V'$ blown up from the isolated fixed points of $\sigma^2$. 
The curve  $\bar W$ is fixed by
$\sigma^2$, while
$\sigma$ fixes
$Z_0$ and two points on $\bar{W}$. 

Since $V'$ is birational to $V$,
it suffices to show that $\pi_1(V')$ is commutative. 
The rational surface $F$ is a blowup of the plane, 
so that $\pi_1 (F - W')$ is generated by
the image of the fundamental group of $L-W'$ where $L$ is the pre-image on $F$ of a line on
${\bb P}^2$. Abusing the terminology we shall still call $L$ a line. Since $V'$ is the
double cover of $F$ along $W'$, $\pi_1 (V' - R'- \cup R_i)$ is a subgroup of
index two in $\pi_1 (F-W')$.
The line $L$ is isomorphic to ${\bb P}^1$,
thus $\pi_1 (L - W')$ is a free group on $n$ generators, where 
$n = \# ({L \cap W'}) -1$.
Generators for $\pi_1 (L-W')$ can be chosen corresponding to small loops 
$\gamma_i$
around the
intersection points of $L$ and $W$, such that the product 
$\gamma_0 \dots \gamma_{n}$ is trivial.
The even products of the $\gamma_i$ generate $\pi_1 (V')$; the relations in
$\pi_1(V')$ are
${\gamma_i}^2 =1$ plus those coming from $\pi_1 (F-W')$.

We will choose the line $L$ on $F$ by considering the fibres of the map 
$\varphi \; : V \to {\bb P}^1$  associated
to the genus two pencil $|3K_V - R|$.  
This is the pullback of the pencil of lines in the plane through $q$,
the order four point of the branch curve $W$.
The member of this pencil through $Q_0$, the isolated fixed point
of $\sigma$,  is decomposed as 
$A_0 + {A_0}'$, the sum of two elliptic $(-1)$ curves which are
interchanged by $\sigma$.  
The curves $A_0$ and ${A_0}'$
meet only at $Q_0$, and
each intersects $R$ in three points. The pre-image of $A_0 + {A_0}'$ on $V'$ is the
$\sigma$-invariant divisor $A_0+A_0'+R_0$; here we use the same notation for the
proper transforms of $A_0$, $A_0'$ on $V'$. 

Set $\ell = A_0 /\sigma^2$; then $\sigma \ell = {A_0}' / \sigma^2$ and $\ell + \sigma \ell
+Z_0$ is a $\sigma-$invariant fibre on $F$.  
Since $\sigma$ switches the fibres of $\varphi$,
this fibration also descends to $\bar F= F/\sigma $.
On $\bar F$, $\ell$ and $\sigma \ell$ become identified.  Set 
$\bar{\ell}= (\ell + \sigma \ell)/\sigma$ and $\bar{W}' = W'/\sigma$.

\proclaim Lemma 5.4. ({\bf [KC]}, Theorem~2)  
Let $Y \to {\bb P}^1$ be a fibration of a quasi-projective variety
$Y$ such that every fibre
 has a smooth component.  
Then for every fibre
$L$ the sequence
$$ \pi_1 (L ) \to \pi_1 (Y) \to \pi_1 ({\bb P}^1 ) =1$$
is exact.

\bigskip

Thus we obtain a surjective homomorphism of groups
$$\pi_1 (\ell + \sigma \ell  - W') \to \pi_1 (F - W').$$ 
\medskip
\proclaim Lemma 5.5.
The map 
$$\pi_1 (\ell - W') \to \pi_1 (F- W')$$
is surjective.

{\sl Proof.}
Consider the diagram
$$
\matrix{
\pi_1(\ell -W') &\to &\pi_1 (F - W') \cr
 \downarrow              &   &\downarrow          \cr
 \pi_1(\bar{\ell} - \bar{W}') &\to &\pi_1 (\bar{F} -\bar{W}') \cr
}$$

By Lemma~5.4, the bottom horizontal map is surjective.   
Since $F-W'$ is an unramified cover of $\bar{F}-\bar{W}'$,
$\pi_1 (F-W') \to \pi_1 (\bar{F} - \bar{W}')$ is injective.
Also the left vertical map is surjective, since 
$\ell -W'$ is isomorphic to $\bar \ell - \bar{W}'$.
Thus any element in $\pi_1 (F-W')$ can be lifted to $\pi_1 (\ell -W')$.
\bigskip

Thus the fundamental group of $F-W'$ is equal to the image of 
$\pi_1(\ell -W')$.
Next consider the commutative diagram
$$
\matrix{
\pi_1(A_0 -R') &\to &\pi_1 (V' - R'- \cup R_i) \cr
 \downarrow              &   &\downarrow \cr
 \pi_1(\ell-W') &\to &\pi_1 (F -W'). \cr
}$$
The vertical arrows are injections onto a subgroup of index 2. By Lemma~5.5 the bottom
horizontal arrow is surjective. Thus the top horizontal arrow is surjective. Finally
consider the commutative diagram
$$
\matrix{
\pi_1(A_0 -R') &\to &\pi_1 (V' - R'- \cup R_i) \cr
 \downarrow              &   &\downarrow \cr
 \pi_1(A_0) &\to &\pi_1 (V') \cr
}$$
corresponding to the natural embeddings of the corresponding varieties. Since
the vertical arrows are surjective, we obtain that the lower horizontal arrow is surjective. 
Thus the fundamental group of $V'$ is the image of the fundamental
group of an elliptic curve, which is commutative.  Therefore
$\pi_1(V) = H_1 (V, {\bb Z}) = \{1\}$, and $V$ is simply connected.

\vglue .5in
\centerline{{\bf 6. Comparison with Barlow surfaces}}
\medskip
Here we prove that our surfaces are not isomorphic to Barlow surfaces. Since any Barlow
surface contains four $(-2)$-curves it suffices to show that $V$ does not
contain any such curves.

\medskip\noindent

\proclaim Lemma 6.1.
({\bf[B]}, Proposition~1.3)
If $S$ is a numerical Campedelli surface (i.e. $p_g = 0, K_{S}^2 = 2$), 
$V$ is a  simply connected surface,
and
$S \rightarrow V$ is a double cover branched along four $(-2)-$curves,
then $S$ has no even torsion.

\medskip\noindent

\proclaim Theorem 6.2. $V$ contains no $(-2)$-curves and hence $V$ is not isomorphic
to a Barlow surface. 

{\sl Proof.} By Lemma 5.1, the set of nodal curves consists 
of four disjoint
curves $C, \sigma(C),$ $ \sigma^2(C),$ and $\sigma^3(C)$. We have
$$C + \sigma C + \sigma^2 C + \sigma^3 C \sim 2(3K_V-R)-2A-2A' = 2L,$$
where $A,A'$ are two disjoint curves of arithmetic genus 1 with 
$A^2 = A'^2 = -1$.   

  We take
the double cover of $V$ branched along the four $(-2)-$curves, $f : Z \rightarrow V$.
We have
$$ L^2 = {{1}\over{4}} \sum_0^3 {\left(\sigma^i C\right)}^2 = -2 
\hbox{ and }
L \cdot K_V =0,$$
thus by Riemann-Roch, the Euler characteristic of $Z$ is
$$\chi_Z = 2 \chi_V +{{L \cdot (L+K_V)}\over{2}} =1.$$
We also have $${K_Z}^2 = 2 (K_V +L)^2 = -2.$$
On $Z$ the branch curves $\sigma^i C$ are disjoint $(-1)-$curves;
blowing these four curves down gives a surface $Z'$ with
${K_{Z'}}^2 =2$.
There are two possibilities for $Z'$: either 
 $p_g=q=1$, or
$p_g =q =0$, in which
case $Z'$ is a numerical Campedelli surface.  We next 
show that neither case can occur, and therefore $V$ cannot 
contain a $(-2)-$curve.

First suppose $p_g = q =1$.  
Since $p_g(Z)= p_g (V) + h^0 (V, {\cal O}_V (K_V+L) ) $ and $p_g (V)=0$,
the space
$H^0 ( {\cal O}_V (K_V+L) )$ is one-dimensional.  We have
$$ (K_V + L)\cdot \sigma^i(C) = L\cdot \sigma^i(C) = {1\over 2}(\sum_0^3 \sigma^i C)\cdot
\sigma^iC = -1.$$
This implies that the unique divisor from $|K_V + L|$ contains all nodal curves
$\sigma^i(C)$. But then $K_V+L-\sum_0^3 \sigma^i(C) = K_V+L-2L = K_V-L$ is effective. 
We have $(K_V-L)^2 = -1$ and $K_V\cdot (K_V-L) = 1$. Thus $K_V-L$ is represented by a curve
$M$ of arithmetic genus $1$. Also, from the proof of Lemma~5.1, each elliptic curve $E_i$ intersects exactly
two of the curves
$\sigma^i(C)$. This implies $(K_V-L)\cdot E_i = 0$. If $E_i$ is an irreducible
component of the unique effective divisor $D\in |K_V-L|$, then all of the $E_i$'s are its
irreducible components. Intersecting with $K_V$ we see that it is impossible. Hence
$M$  is mapped isomorphically onto a curve $\bar M$
in the quintic model. Since, again by Lemma 5.1, $\pi^*(H)\cdot \sigma^i(C) = 2$ we
have 
$\pi^*(H)\cdot (K_V-L) = 5-4 = 1$. Thus $\bar M$ is a line and cannot be of arithmetic genus
1.   Therefore $Z'$ cannot be a surface
with $p_g=q=1$.

Next consider the case $K_{Z'}^2 =2$, $p_g=q=0$.  Then the minimal
model $Z'$ is a numerical Campedelli surface.  
By Corollary~3.6, the pencil $\left| 3K-R \right|$ has five exceptional
fibres which contribute to the Picard group $8$ linearly
independent divisors contained in fibres, namely
the four $(-2)$-curves $C_1, \dots, C_4$ contained in two fibres,
three disjoint $(-1)$-elliptic curves $A_1, A_2, A_3$ each
contained in the other three exceptional fibres, plus the 
general fibre. It follows from the proof of Theorem 3.7 that each elliptic
curve $E_i$ intersects $A_i$ or its complementary component $A_i'$ with multiplicity
1. Also we know that each $A_i$ is $\sigma^2$-invariant. Thus we can choose the $A_i$
in such a way that the elliptic curves $E_1$ and $E_2 =\sigma^2(E_1)$ intersect each
$A_i$ with multiplicity $1$.  The Picard number of
$V$ is nine. This implies that $E_1$ and $E_2$ (being 3-sections) differ by a
${\bb Q}$-divisor contained in fibres. Write 
$$E_1 =
E_2+n_1A_1+n_2A_2+n_3A_3+m_1C_1+m_2C_2+m_3C_3+m_4C_4+aF,$$
 where 
$F$ is a general fibre of $|3K-R|$.  We also may assume that 
$$C_2 = \sigma^2(C_1), \,
C_4=\sigma^2(C_3), \, E_1\cdot C_1 = E_1\cdot C_3 = 1, \,
E_2\cdot C_2 = E_2\cdot
C_4 = 1 .$$
Intersecting both sides with the $A_i$ and $C_i$ gives
$$E_1 = E_2 -\frac{1}{2} C_1 -\frac{1}{2} C_2 +\frac{1}{2} C_3 +\frac{1}{2} C_4 ,$$
thus $2 E_1 +C_1 +C_2 \sim 2 E_2 +C_3 +C_4$ on $V$.  On the double cover $S$ of $V$,
we have $2 E_1' +2C_1' +2C_2' \sim 2 E_2' +2C_3' +2C_4'$ where
we write $E_i', C_i'$ for the pullback of $E_i, C_i$ to the cover.  By 
Lemma~6.1, $S$ has no $2$-torsion, thus $E_1' +C_1' +C_2' \sim E_2' +C_3'
+C_4'$, and $\left| E_1'+C_1' +C_2' \right|$ is a genus two pencil on $Z$.  

Now let $A_1$ be an elliptic component of another reducible fibre of
$\left| 3K_V -R \right|$; since $A_1$ does not meet the $(-2)$-curves $C_i$,
it splits on the cover, and gives a section of the pencil $\left| E_1' +C_1'
+C_2'\right|$  on $Z$.  This implies that $A_1$ is an elliptic 
section of a rational pencil, a contradiction.
Thus we cannot obtain such a covering of $V$, and again we see that there
cannot be four $(-2)-$curves on $V$.

\vglue .3in

\centerline{{\bf References}}

\medskip
\noindent
{\bf [B]} Barlow, R., {\it A simply connected surface of general type with 
$p_g = 0$}, Invent. Math., 79 (1985), 293-301. 

\medskip
\noindent
{\bf [CG]} Craighero, P., Gattazzo, R.,
{\it Quintic surfaces of ${\bb P}^3$ having a non singular
model with $q=p_g=0$, $P_2 \neq 0$,}
Rend. Sem. Mat. Univ. Padova, 91 (1994), 187-198.

\medskip
\noindent
{\bf [CL]} Catanese, F., LeBrun, C.,
{\it On the scalar curvature of Einstein manifolds,}
preprint.

\medskip
\noindent
{\bf [KC] }
Keum, J.H., Chung, I.J.,
{\it On fundamental groups of fibred complex manifolds,}  Mich. Math. Journ., to appear.

\medskip
\noindent
{\bf [M] }
Miyaoka, Y.,
{\it Tricanonical maps of numerical Godeaux surfaces,}
Invent. math., 34 (1976) , 99-111.

\medskip
\noindent
{\bf [OP]}
Oort, F., Peters, C.,
{\it A Campedelli surface with torsion group ${\bb Z}/2 {\bb Z}$,}
Indag. Math., 43 (1981), 399-407.

\medskip
\noindent
{\bf [R]} Reid, M., {\it Campedelli versus Godeaux,} in ``Problems in the theory of surfaces and their classification (Cortona, 1988)'', pp. 309-365, Academic
Press, 1991.

\medskip
\noindent
{\bf [St] }
Stagnaro, E.,
{\it On Campedelli branch loci,}
Dip. di Met. e Mod. Matematici, Universita di Padova,  53, (1996).

\medskip
\noindent
{\bf [W] } 
Werner, C.,  
{\it  A surface of general type with $p_g=q=0$, $K^2=1$,}
Manuscripta math., 84 (1994), 327-341.

\vglue .4in
{\sevenrm DEPARTMENT OF MATHEMATICS, UNIVERSITY OF MICHIGAN, 
ANN ARBOR, MI 48109}

{\it E-mail address}: idolga@math.lsa.umich.edu, cwerner@math.lsa.umich.edu
\bye